\documentclass[12pt]{article}

\usepackage{amsmath,amsthm}
\usepackage{amssymb}
\usepackage{color}
\usepackage{xspace}
\usepackage[colorlinks=true,
linkcolor=green,
filecolor=brown,
citecolor=green]{hyperref}

\def\a{\ensuremath{\mathcal A}\xspace}

\def\gl{ground level\xspace}

\begin{document}
 
\newtheorem{theorem}{Theorem}
\newtheorem{prop}[theorem]{Proposition}
\newtheorem{cor}[theorem]{Corollary}
\mbox{ }
\vspace*{-25mm}
\begin{center}
{\Large
The 136th Manifestation of $C_{n}$                          \\ 
}

\vspace{10mm}
DAVID CALLAN  \\
Department of Statistics  \\
\vspace*{-2mm}
University of Wisconsin-Madison  \\
\vspace*{-2mm}
1300 University Ave  \\
\vspace*{-2mm}
Madison, WI \ 53706-1532  \\
{\bf callan@stat.wisc.edu}  \\
\vspace{3mm}

October 31, 2005
\end{center}

\begin{abstract}
    We show bijectively that the Catalan number $C_{n}$ counts Dyck 
    $(n+1)$-paths in which the terminal descent is of even length and all 
    other descents to ground level (if any) are of odd length.
\end{abstract}

\vspace*{2mm}

Richard Stanley's 
\htmladdnormallink{inventory }{http://www-math.mit.edu/~rstan/ec/} of 
combinatorial interpretations of the Catalan 
number
\htmladdnormallink{$C_{n}$}{http://www.research.att.com:80/cgi-bin/access.cgi/as/njas/sequences/eisA.cgi?Anum=A000108}
$\,=\frac{1}{n+1}\binom{2n}{n}$ currently stands at 135 items. Here 
is one more.

\begin{theorem}
    Let $\a_{n}$ denote the set of Dyck $n$-paths for which the terminal descent is of even length 
    and all other descents to ground level $($if any$)$ are of odd 
    length. Then $\vert \a_{n}\vert = C_{n-1}$ for $n\ge 2$.
\end{theorem}
\vspace*{-4mm}
This result is a counterpart to item ($j$) in Stanley's inventory, 
which says that $C_{n-1}$ also counts Dyck $n$-paths for which \emph{all} 
descents to \gl are of odd length.

A Dyck $n$-path is a lattice path of $n$ upsteps $U$ and $n$ downsteps $D$ 
that never dips below \emph{ground level}, the horizontal line joining its 
start and end points. The number of Dyck $n$-paths is well known to be $C_{n}$. The \emph{size}, also called the 
\emph{semilength}, of a Dyck $n$-path is $n$. A \emph{return} is a downstep that 
returns the path to \gl. 
A \emph{descent} is a maximal sequence of contiguous 
downsteps. A \emph{peak} is an occurrence of $UD$. A \emph{low} peak (resp. low $UDU$) is one that 
starts at \gl. A low peak is also called a \emph{hill} and a low $UDU$ an \emph{early hill}. 
Note that a path free of early hills is either hill-free 
or has just one hill at the very end. Hill-free Dyck paths and Dyck 
paths with an even-length terminal descent are both counted \cite{fine} by the 
Fine numbers,
\htmladdnormallink{A000957}{http://www.research.att.com:80/cgi-bin/access.cgi/as/njas/sequences/eisA.cgi?Anum=A000957}
in OEIS. Early-hill-free Dyck paths are counted \cite{udu} by
\htmladdnormallink{A000958}{http://www.research.att.com:80/cgi-bin/access.cgi/as/njas/sequences/eisA.cgi?Anum=A000958}
.   

We prove the following refinement of $\vert \a_{n}\vert = C_{n-1}$.

\newpage 
\mbox{ }
\vspace*{-15mm}
\begin{theorem}
    For $n\ge 2$ and $k\ge 1$, the paths in $\a_{n}$ with $k$ returns correspond 
    bijectively to Dyck $(n-1)$-paths that contain $k-1$ early hills.
\end{theorem}
The proof relies on the following bijections.

\begin{prop}
    There exists a bijection from Dyck $n$-paths with terminal descent 
    of even $($resp. odd\,$)$ length to hill-free $($resp. 
    early-hill-free$)$ 
    Dyck $n$-paths.
\end{prop}

\noindent \textbf{Proof}\quad The ``\emph{DUtoDXD}'' bijection of 
\cite[\S 4]{some} establishes the even-length terminal descent $\rightarrow$ 
hill-free part.
For the odd-length terminal descent $\rightarrow$ early-hill-free 
part, split the first set 
of paths into $A$: those with only one return, and $B$: those with 2 
or more returns. The interior (drop first and last steps) of a path 
in $A$ has terminal descent of even length and so corresponds to a 
hill-free Dyck $(n-1)$-path by the previous part. Append $UD$ to get 
a bijection from $A$ to
the early-hill-free Dyck $n$-paths that end $UD$. A path in $B$ can be  
\linebreak\\[-4mm]
written (uniquely) as $PUQD=
P\ \diagup\ \raisebox{4mm}{$Q$}\ \diagdown\ $
where $P,Q$ are nonempty Dyck paths and \linebreak\\[-4mm] 
$Q$ has terminal descent of even length. Map to
$\:\diagup\ \raisebox{3mm}{$P$}\ 
\diagdown\ Q'\:$,
 where $Q'$ is the hill-free path corresponding to $Q$. This gives 
a bijection from $B$ to the early-hill-free Dyck $n$-paths that do not 
end $UD$. \qed

\noindent \textbf{Proof of Theorem 2}\quad Given a path in $\a_{n}$ 
with $k$ returns, use the path's 
returns to \linebreak\\[-4mm]
write it (uniquely) as 
$
\ \diagup\ \raisebox{4mm}{$P_{1}$}\ \diagdown\
 \diagup\ \raisebox{4mm}{$P_{2}$}\ \diagdown\quad \ldots \quad
 \diagup\ \raisebox{4mm}{$P_{k-1}$}\ \diagdown\
 \diagup\ \raisebox{4mm}{$P_{k}$}\ \diagdown\
$ 
 where $P_{1},P_{2},\ldots,P_{k-1}$ are Dyck paths, all with 
terminal descent of even length (possibly 0), and $P_{k}$ is  
a \linebreak
Dyck path with terminal descent 
of odd length. Using Prop.\,3, map the path to \linebreak
$P_{1}'\ \diagup\,\diagdown\
 P_{2}'\ \diagup\,\diagdown\quad  \ldots \quad 
 \diagup\,\diagdown\
 P_{k-1}'\ \diagup\,\diagdown\
P_{k}'\,$, 
where $P_{i}'$ is hill-free for $1\le i 
\le k-1$ and $P_{k}'$ is nonempty early-hill-free. The resulting Dyck path 
has one fewer $U$ and $D$ than the original and contains $k-1$ early 
hills, and 
Theorem 2 follows. \qed

These results can be used to explain the distribution of the 
statistic ``\#\ even-length descents to \gl'' on Dyck paths. First, let $T(n,k)$ denote the 
number of Dyck $n$-paths with $k$ returns; $\big(T(n,k)\big)_{0\le k 
\le n}$ forms the Catalan triangle, 
\htmladdnormallink{A106566}{http://www.research.att.com:80/cgi-bin/access.cgi/as/njas/sequences/eisA.cgi?Anum=A106566}
in OEIS. 
\begin{cor}[\cite{sun2}]
    The number of Dyck $n$-paths with $k$ even-length descents to \gl 
    is $T(n,2k)+T(n,2k+1)$.
\end{cor}
\textbf{Proof}\quad Again calling on the ``\emph{DUtoDXD}'' bijection of 
\cite[\S 4]{some}, it sends Dyck $n$-paths all of whose returns to \gl 
have odd length to Dyck $n$-paths that start $UD$ and thence 
(transfer this $D$ to the end of the path) to Dyck $n$-paths with 
exactly 1 return. This establishes the case $k=0$. For $k \ge 1$, 
split the paths into $A$: those for which the terminal descent has 
even length, and $B$: the rest. A path in $A$ splits, via its 
even-length descents to \gl, into $k$ Dyck paths to each of which Theorem 1 
applies. The result is a $k$-list of \emph{nonempty} Dyck paths of 
total size $n-k$. Since nonempty Dyck paths \\[2mm] correspond to 2-return 
Dyck paths of size 1 unit larger \big(\:$\diagup\ \raisebox{3mm}{$P$}\ 
\diagdown\ Q \to \diagup\ \raisebox{3mm}{$P$}\ 
\diagdown\ \diagup\ \raisebox{3mm}{$Q$}\ 
\diagdown$\:\big), we get a bijection from $A$ to Dyck $n$-paths with $2k$ 
returns. There is a similar bijection from $B$ to Dyck $n$-paths with 
$2k+1$ returns. \qed

\end{document}